# Cycles and Patterns in the Sieve of Eratosthenes—Part 2, Potential Twin Primes

George F. Grob

June 14, 2021

This paper analyzes the emergence and distribution of potential twin primes, i.e., pairs of integers, $a - 1$ and $a + 1$, that are both relatively prime to the first $n$ primes or to a given set $M$ of primes, and which are the breeding grounds of true twin primes. It describes cyclical patterns of their location across the number line and provides a formula for counting them.

**1. Introduction**

In the third century B.C.E. Eratosthenes of Cyrene (276–194 BC) invented his famous sieve to identify prime numbers. Even today it remains an irresistible and practical starting point for analyzing them. Almost all mathematicians who have explored it have noticed that after each iteration there remain pairs of integers separated by 2. Some such pairs survive as twin primes; others do not, as one or both members are divided by a prime in subsequent iterations of the sieve. However, there are always an infinite number of surviving twin pairs, albeit increasingly rare. This has fueled the as yet unproven twin prime conjecture, that there are an infinite number of twin primes.

The challenge of proving the twin prime conjecture has been the fixation and bane of mathematicians for over two thousand years. Remarkably, though, less attention has been paid to Eratosthenes's Sieve itself, particularly the "second" part, the survivors left over after each iteration of it. That is our primary focus here, particularly the relentlessly surviving twin pairs.

**A. Primes, *N*-Primes, *M*-Primes, and Twins.** Primes are integers other than the integer 1 that are not divisible by any integer other than themselves and 1. Analogously, "$n$-primes" are integers not divisible by any of the first $n$ primes, and "$M$-primes" are integers not divisible by any member of an arbitrarily given set, $M$, of primes. In context, either of these may be referred to as "potential primes."

*Twin primes, twin n-primes,* and *twin M-primes* are, respectively, pairs of primes, $n$-primes, or $M$-primes whose difference is 2. *These twins are the focus of this paper*.

A *twin prime center* is an integer, x, such that x – 1 and x + 1 are both primes. Analogously, a *twin n-prime or M-prime center* is an integer, x, such that x – 1 and x + 1 are both n-primes or both M-primes.

This paper follows up on a broader analysis of potential primes by this author and M. Schmitt, *Cycles and Patterns in the Sieve of Eratosthenes.*[i] (available at https://arxiv.org/abs/1905.03117). The key theorems of that paper are listed in the appendix. Among them, the following two, which focus on "twins," are the more specific starting points for this paper.

- Theorem 4. Twin n-primes are pairs consisting of two $n$-primes, $(a - 1)$ and $(a +1)$. The number of such pairs $\leq \Pi p_i = \Pi(p_i - 2)$, $1 \leq i \leq n$ ($p_i \neq 2$).
- Theorem 4a. Twin M-primes are pairs consisting of two M-primes, $(a – 1)$ and $(a +1)$. The number of such pairs $\leq \Pi p_i = \Pi(p_i - 2)$, ($p_i \, \varepsilon \, M$).

**B. Questions to be Analyzed.** Following are the questions and key concepts around which we have organized our findings.

1) *Sieve of Twins.* Is there a sieve, analogous to the Sieve of Eratosthenes, that incrementally produces lists of twin $n$-primes or twin $M$-primes?
2) *Legendre-like Function for Counting Twins.* Is there an algebraic formula for counting twin $n$-primes or twin $M$-primes analogous to Legendre's well-known formula for counting primes[ii]?
3) *Meissel/Lehmer-like Step-Wise, Inductive Twin Counting Function.* Is there an inductive, step-wise function for counting twins that is analogous to Meissel/Lehmer's well-known inductive formula for calculating the number of integers $\leq$ x that are relatively prime to the first $n$ primes based on the number of integers relatively prime to the first ($n$ - 1) primes[iii]?
4) *Euler-like ϕ-function for Twins.* Is there a formula analogous to Euler's ϕ-function[iv], that for any integer x returns the number of pairs of integers, $(a - 1)$ and $(a + 1)$, that are $\leq$ x with both relatively prime to x?

## 2. Background

**A. Deeper Look at the Sieve of Eratosthenes.** Eratosthenes's Sieve methodically and relentlessly produces a list of prime numbers, and does so in the order in which they appear. Each new prime is the first survivor of the of the previous iteration of his sieve. Given any number, x, subject to the limits of time and computing capacity, Eratosthenes's Sieve will ultimately reach and single out every prime number that is $\leq$ x, as well as each pair of twin primes.

As widely known and intuitively clear as it is, it is helpful to step back for a wider view in order to set the stage for a "Twin Sieve" that identifies twin.

With Eratosthenes's Sieve, each phase, or what we may more appropriately call "cycle," starts with the first $n$ primes . . . i.e., 2, 3, 7, 11, . . . $n$. After the execution of the "$n^{th}$" phase, the sieve contains two parts: 1) a list of the first $n$ primes, and 2) an infinitely long residue of integers, all of them relatively prime to the first $n$ primes. The first remnant is the $(n + 1)^{st}$ prime. This is transferred to the first part and the process of discarding its multiples in the second part is repeated, and so on, as many times as one wishes. Doing so produces a list of prime numbers, and does so in the order that they appear along the number line.

**B. Identifying and Counting Twins.** Eratosthenese's Sieve does not, per se, produce lists of *twin* primes. It ultimately reaches any such pairs; but does not identify or highlight them as such. Additional steps are needed for that. Furthermore, with respect to twins, Eratosthenese's Sieve is inherently inefficient. The remnant of each cycle mixes potential *twin* primes with *all* potential primes.

Therefore, we introduce and describe here an additional sieve, a "Sieve of Twins," modeled on Eratosthenes's, but focused on what we call "twin $n$-primes" i.e. integers $a – 1$ and $a + 1$ that are potentially twin primes, both being relatively prime to the first *n* primes. Analogously, our Sieve of Twins can be applied to *M*-primes.

Just as each cycle of Eratosthenes's Sieve discards integers that cannot possibly be primes, our Sieve of Twins sequentially discards endless lists of integers that cannot possibly be twin prime centers. In both cases, Eratosthenes's Sieve and our Sieve of Twins, with each succeeding "*n*" there are always remnants. Those remnants become thinner and thinner, but remain infinitely long, containing, respectively, potential primes and potential twin prime centers.

**3. Our Sieve of Twins**

Broadly, our Sieve of Twins mimics Eratosthenes's Sieve of Primes. However, instead of sequentially identifying primes and then discarding their multiples, it starts with a set of primes (e.g., the first *n* primes, or any set *M* of primes) and then incrementally discards all integers, x, such that either (x – 1), (x + 1), or both, are divisible by one or more of the primes in the initial set. What remains is an infinite set of integers, x, that are the centers of twin integers, (x – 1) and (x + 1), both of which are relatively prime to our selected set of primes.

For any single prime, $p_j$, all integers, x, are centers of twin $\{p_j\}$-primes, (x - 1) and (x + 1), except when x = ($p_j$k – 1) or ($p_j$k + 1), for some integer, k.

If x takes the form ($p_j$k – 1), then x + 1 = $p_i$·k and is therefore *not* a $\{p_j\}$-prime. Similarly, if x takes the form ($p_j$k + 1), then x – 1 = $p_i$·k, and it is therefore *not* a $\{p_j\}$-prime. In either or both cases, the pair, ($p_j$k – 1), ($p_j$k + 1), is *not* a pair of twin n-primes. This phenomenon is illustrated in *Table 1, Centers of Twin {7}-Primes*, for twins relatively prime to 7.

**Table 1. Centers of Twin {7}-Primes**

| Multiples of 7 | {7}-Primes | Twin {7}-Primes | Twin {7}-Prime Centers | Not Twin {7}-Prime Centers | Number of Twin {7}-Prime Centers |
|---|---|---|---|---|---|
| 0 | | | | | |
| | 1 | | | 1 | |
| | 2 | 1, 3 | 2 | | 1 |
| | 3 | 2, 4 | 3 | | 2 |
| | 4 | 3, 5 | 4 | | 3 |
| | 5 | 4, 6 | 5 | | 4 |
| | 6 | | | 6 | |
| 7 | | 6, 8 | 7 | | 5 |
| | 8 | | | 8 | |
| | 9 | 8, 10 | 9 | | 6 |
| | 10 | 9, 11 | 10 | | 7 |
| | 11 | 10, 12 | 11 | | 8 |
| | 12 | 11, 13 | 12 | | 9 |
| | 13 | | | 13 | |
| 14 | | 13, 15 | 14 | | 10 |
| | 15 | | | 15 | |
| | 16 | 15, 17 | 16 | | 11 |
| | 17 | 16, 18 | 17 | | 12 |
| | 18 | 17, 19 | 18 | | 13 |
| | 19 | 18, 20 | 19 | | 14 |
| | 20 | | | 20 | |
| 21 | | 20, 22 | 21 | | 15 |
| | 22 | | | 22 | |
| | 23 | 22, 24 | 23 | | 16 |
| etc. | | | | | |

The table also illustrates that, in accordance with Theorem 1 in Appendix 1, for a single prime, $p$, (in this case 7), the number of twin prime centers is cyclical, repeating every $p$ (here, 7) integers.

Tables 2 through 6 illustrate the sequential "killing off" of potential twin primes by primes 2, 3, 5, and 7. Again, in accordance with Theorem 1, the cycle for the first prime, 2, repeats every 2 integers; the combined effect of the first two primes, 2 and 3, repeats every 6 times; the combined effect of the first three primes, 2, 3, and 5, repeats every 30 times; and the combined effect of the first four primes, 2, 3, 5, and 7, repeats every 210 times. Thus, it is convenient and revealing to examine the first 210 integers to illustrate the successive impacts of each of the first four primes with respect to our Sieve of Twins (hereafter "Twin Sieve"). In turn, what happens in these first four cycles foreshadows the broader patterns that emerge as a result of applications of the Twin Sieve for subsequent primes, such as 11, 13, 17, etc. (Again, see Theorem 1 in Appendix 1.)

We begin, then, by illustrating what happens within the first 210 integers, encompassing the first cycle for each of the first 4 $n$–primes, as well as their repetitions for the first three.

**Table 2. The First Cycle of the First 4 $n$-Primes (2, 3, 5, 7)**

| 1 | 2 | 3 | 4 | 5 | 6 | 7 | 8 | 9 | 10 |
|---|---|---|---|---|---|---|---|---|----|
| 11 | 12 | 13 | 14 | 15 | 16 | 17 | 18 | 19 | 20 |
| 21 | 22 | 23 | 24 | 25 | 26 | 27 | 28 | 29 | 30 |
| 31 | 32 | 33 | 34 | 35 | 36 | 37 | 38 | 39 | 40 |
| | | | | etc. | | | | | |
| 191 | 192 | 193 | 194 | 105 | 196 | 197 | 198 | 199 | 200 |
| 201 | 202 | 203 | 204 | 205 | 206 | 207 | 208 | 209 | 210 |

Twin "1-primes" are pairs of integers relatively prime to the first prime, 2. Thus, the first cycle of the Twin Sieve drops integers of the form $(2k – 1)$ and $(2k + 1)$ from Table 2. (In this special case of 2, these two cycles overlap for all $k \geq 2$.) The remaining integers, in this case all the even numbers, are the centers of twin 1-prime pairs, such as (1, 3); (3, 5); (5, 7), etc., as shown in Tables 3A and 3B.

**Table 3. Twin 1-Primes**

| A. Twin 1-Prime Centers | | | | | B. Twin 1-Prime Pairs | | | | |
|---|---|---|---|---|---|---|---|---|---|
| 2 | 4 | 6 | 8 | 10 | 1, 3 | 3, 5 | 5, 7 | 7, 9 | 9, 11 |
| 12 | 14 | 16 | 18 | 20 | 11, 13 | 13, 15 | 15, 17 | 17, 19 | 19, 21 |
| 22 | 24 | 26 | 28 | 30 | 21, 23 | 23, 25 | 25, 27 | 27, 29 | 29, 31 |
| 32 | 34 | 36 | 38 | 40 | 31,33 | 33, 35 | 35, 37 | 37, 39 | 39, 41 |
| | | etc. | | | | | etc. | | |
| 192 | 194 | 196 | 198 | 200 | 191, 193 | 193, 195 | 195, 197 | 197, 199 | 199, 201 |
| 202 | 204 | 206 | 208 | 210 | 201, 203 | 203, 205 | 205, 207 | 207, 209 | 209, 211 |

Twin "2-primes" are pairs of integers relatively prime to the first two primes, 2 and 3. Thus, the second iteration of the Twin Sieve drops integers that remain from the first culling but which are of the form $(3k – 1)$ or $(3k + 1)$. The remnants of this second culling are shown in Tables 4A and 4B.

**Table 4. Twin 2-Primes**

| A. Twin 2-Prime Centers | | | | | B. Twin 2-Prime Pairs | | | | |
|---|---|---|---|---|---|---|---|---|---|
| 6 | 12 | 18 | 24 | 30 | 5, 7 | 11, 13 | 17, 19 | 23, 25 | 29, 31 |
| 36 | 42 | 48 | 54 | 60 | 35, 37 | 41, 43 | 47, 49 | 53, 55 | 59, 61 |
| 66 | 72 | 78 | 84 | 90 | 65, 67 | 71, 73 | 77, 79 | 83, 85 | 89, 91 |
| 96 | 102 | 108 | 114 | 120 | 95, 97 | 101, 103 | 107, 109 | 113, 115 | 119, 121 |
| 126 | 132 | 138 | 144 | 150 | 125, 127 | 131, 133 | 138, 139 | 143, 145 | 149, 151 |
| 156 | 162 | 168 | 174 | 180 | 155, 157 | 161, 163 | 167, 169 | 173, 175 | 179, 181 |
| 186 | 192 | 198 | 204 | 210 | 185. 187 | 191, 193 | 187, 189 | 203, 205 | 209, 211 |

Twin "3-primes" are pairs of integers relatively prime to the first three primes (2, 3, and 5). Thus, in the same vein, the third culling of the Twin Sieve drops integers that remain from the second iteration but are of the form $(5k – 1)$ or $(5k + 1)$. The remnants of this third culling are shown in Tables 5A and 5B.

**Table 5. Twin 3-Primes, Centers and Pairs**

| A. Twin 3-Prime Centers | | | B. Twin 3-Prime Pairs | | |
|---|---|---|---|---|---|
| 12 | 18 | 30 | 11, 13 | 17, 19 | 29, 31 |
| 42 | 48 | 60 | 41, 43 | 47, 49 | 59, 61 |
| 72 | 78 | 90 | 71, 73 | 77- 79 | 89, 91 |
| 102 | 108 | 120 | 101, 103 | 107, 109 | 119, 121 |
| 132 | 138 | 150 | 131, 133 | 137. 139 | 149, 151 |
| 162 | 168 | 180 | 161, 163 | 167, 169 | 179, 181 |
| 192 | 198 | 210 | 191, 193 | 197, 199 | 209, 211 |

Twin "4–primes" are pairs of integers relatively prime to the first four primes (2, 3, 5, and 7). Following the same pattern as above, the fourth culling of the sieve drops from the remnants of the third those integers that are of the form $(7k - 1)$ and $(7k + 1)$, shown in Tables 6A and 6B.

**Table 6. Twin 4-Primes, Centers and Pairs**

| A. Twin 4-Prime Centers | | | B. Twin 4-Prime Pairs | | |
|---|---|---|---|---|---|
| 12 | 18 | 30 | 11, 13 | 17, 19 | 29, 31 |
| 42 | 60 | 72 | 41, 43 | 59, 61 | 71, 73 |
| 102 | 108 | 138 | 101, 103 | 107, 109 | 137, 139 |
| 150 | 168 | 180 | 149, 151 | 167, 169 | 179, 181 |
| 192 | 198 | 210 | 191, 193 | 197, 199 | 209, 211 |

**E. Comparing the Sieves of Eratosthenes and Twins.** With this explanation and examples in mind, it may be helpful to compare the Sieve of Eratosthenese's and our Sieve of Twins. They have much in common, but there are some differences. Both of them are cyclical and each has two parts, as highlighted in the following table.

**Table 7. The Sieves of Eratosthenes's and Twins Compared**

|  | Part 1, Primes | Part 2, Survivors |
|---|---|---|
| **Sieve of Eratosthenes** | First $n$ primes | N-primes |
| **Sieve of Twins** | First $n$ primes | Twin N-prime centers |

**Eratosthenes's Sieve** sequentially and relentlessly produces a list of primes, and does so in numerical order, from the first, second, and to the $n^{th}$. They are listed sequentially in Part 1 of the sieve.

What is left over after each cycle, Part 2, is a list of what we call $n$-primes—i.e., integers that are relatively prime to the first $n$ primes (that are listed in Part 1 of the sieve). The n-primes appear as subsets within infinitely repetitive, symmetric cycles of integers that each form modular rings with multiplicative and additive algebraic functions. Within these rings are also other cycles of integers of uniform length and containing identical numbers of integers relatively prime to the first *n* primes[v].

**The Sieve of Twins** also lists the first *n* primes in Part 1. But Part 2 of this sieve consists of the *centers* of twin *n*-prime pairs. These centers are also symmetrically located within and are subsets of the same recurring modular rings of the n-primes. However, they do not of themselves form groups or rings. Nevertheless, they are of special interest, in that for any integer, *n*, the centers of *all* twin primes greater than *n* are subsets of the centers of twin *n*-prime pairs.

It is tempting, but misleading (and disappointing), to trust that, like Eratosthenes's Sieve, the first survivor of each cycle of the Sieve of Twins will be the center of twin primes. It almost certainly is. But we cannot prove it. It remains what we may wish to call a "potential twin prime center."

On the other hand, just as Eratosthenes's Sieve produces the "next" prime, the Sieve of Twins can easily and sequentially produce a list of the "next" twin primes—or more precisely, a list of the centers of twin prime pairs. The process for doing so is as follows.

Every twin *n*-prime center, *a*, would be is the center of true twin primes, (*a* - 1) and (*a* + 1), if neither of them is divisible by a prime. Since, by definition, neither one can be divided by any of the first *n* primes, the smallest prime that could divide either one would be the $(n + 1)^{st}$ prime, $p_{n+1}$, and the smallest multiple of $p_{n+1}$ that could do so would be $(p_{n+1})^2$. Thus, as long as $a < ((p_{n+1})^2 - 2)$, then (*a* - 1) and (*a* + 1) will both be smaller than $(p_{n+1})^2$, and therefore both will be primes, and hence will be pairs of true twin primes.

This serial derivation of true twin primes from twin *n*-primes is illustrated in *Table 8, Twin n-Primes That Are Also Twin Primes*.

**Table 8. Twin *n*-Prime Pairs That Are Also Twin Primes**
First number in each cell = twin prime center;
Numbers in parentheses = twin prime pairs

| *n* = 1<br>{ 2 } | *n* = 2<br>{ 2, 3 } | *n* = 3<br>{ 2, 3, 5 } | *n* = 4<br>{ 2, 3, 5, 7 } |
|---|---|---|---|
| 4  (3, 5) | | | |
| 6  (5, 7) | 6  (5, 7) | | |
| $p_{n+1} = 3$<br>$(p_{n+1})^2 – 2 = 9$ | 12  (11, 13) | 12  (11, 13) | 12  (11, 13) |
| | 18  (17, 19) | 18  (17, 19) | 18  (17, 19) |
| | $p_{n+1} = 5$<br>$(p_{n+1})^2 – 2 = 23$ | 30  (29, 31) | 30  (29, 31) |
| | | 42  (41, 43) | 42  (41, 43) |
| | | $p_{n+1} = 7$<br>$(p_{n+1})^2 – 2 = 47$ | 60  (59, 61) |
| | | | 70  (71, 73) |
| | | | 102  (101, 103) |
| | | | 108  (107, 109) |
| | | | $p_{n+1} = 11$<br>$(p_{n+1})^2 – 2 = 119$ |

As with the Sieve of Eratosthenes, our "Twin Sieve" can be repeated for each succeeding prime. For example, for the first *n* primes (2, 3, 5, 7, . . . n), it is natural to apply the sieves first to 2, then 3, 5, 7, etc. However, for any given "*n*" we obtain the same results no matter what the order of application. For example, for the first four primes, we get the same results if we first cull out centers of twins relatively prime to, say, 5, then to 3, 7, and 2, etc., in fact, in any order.

In addition, the Twin Sieve can produce lists of twin pairs of any set, *M*, of primes, resulting in endlessly repeating cycles of twins that are relatively prime to the primes in *M*.

By Theorem 2, the centers of such twin *n*-primes or *M*-primes are located symmetrically within repetitive cycles of integers of length $\Pi(p_i)$, $1 \leq i \leq n$, or $p_i \in M$. And, by Theorem 4, within each such cycle, there are $\Pi(p_i - 2)$, $1 \leq i \leq n$, or $p_i \in M$, $p_i \neq 2$, twin pairs. For example, Theorem 4 predicts that there will be $(3-2)(5-2)(7-2) = 1 \cdot 3 \cdot 5 = 15$ centers of twin 4-primes $\leq 2 \cdot 3 \cdot 5 \cdot 7 = 210$, the first such cycle, as seen in Table 68. And, again, by Theorem 1, these result hold true for each subsequent cycle of 210 integers.

**4. A Legendre-like Function for Counting Twins**

**A. Legendre's Formula.** We will adapt Adrien-Marie Legendre's (1752 – 1833) formula for counting primes in order to construct a new formula for counting twin primes, twin *n*-primes, and twin *M*-primes.

Legendre's formula for counting primes is:

(1) $f_n(x) = [x] - \Sigma [x/p_i] + \Sigma [x/p_i \cdot p_j] - \Sigma [x/p_i \cdot p_j \cdot p_k] + \cdots \Sigma [x/p_1 \cdot p_2 \cdot p_3 \cdots p_n]$,
the sums being taken over each one, then over all the doubles, triples, etc. of the first n primes.

For any integer *n* and any real number x, the formula returns the number of integers $\leq x$ that are relatively prime to the first *n* primes. (It is worth noting that if $x < (p_n)^2$, then $f_n(x) - 1$ = the number of primes $p$, $p_n < p < (p_n)^2$. Thus, incrementally increasing the value of *n* and adding that value to the value of Legendre's formula, it can theoretically be used to count the number of all primes of any value for $x < (p_n)^2$.) That is,

(2) $\pi(x) = n + f_n(x) - 1$ for all $x < (p_n)^2$

We may generalize this formula to integers relatively prime to all primes within a given set *M* of primes, i.e.

(3) $f_M(x) = [x] - \Sigma[x/p_i] + \Sigma[x/p_i \cdot p_j] - \Sigma[x/p_i \cdot p_j \cdot p_k] + \cdots \Sigma[x/p_i \cdot p_j \cdot p_k \cdots p_z]$,
the sums taken over each one, then all the doubles, triples, etc. of the primes $\in M$.

While the proof of Legendre's theorem is well known, it is worthwhile to recapitulate it here in order to set the stage for an analogous theorem for twin primes, n-primes, or M-primes that we will introduce below. The proof is based on the "inclusion-exclusion principle."[vi]

By definition, the integer component of any real number x is [x]. Consider any prime, $p_i$. Starting with 0, every $(p_i)^{th}$ integer is divisible by $p_i$. All other integers and every real number have a remainder when divided by $p_i$. Thus, the integral component of any real number x divided by $p_i$ is $[x/p_i]$. So, by definition, the answer to the question, "How many integers $\leq x$ are divisible by $p_i$?" is the integral part of x less the integral component of $x/p_i$, or $[x] - [x/p_i]$.[vii]

With this in mind, consider a set *M* of primes, $p_1, p_2, p_3 \cdots p_k$. Then for a given real number x, the number of integers ≤ x that are relatively prime to $p_1$ is $[x] - [x/p_1]$; the number relatively prime to $p_2$ is $[x] - [x/p_2]$, and so forth, with the number of integers relatively prime to the i[th] primes being $[x] - [x/p_i]$. On the surface then, it would seem that the number of integers ≤ x that are relatively prime to both $p_1$ and $p_2$ would be $[x] - [x/p_1] - [x/p_2]$. However, that would not be accurate, because any multiple of both $p_1$ and $p_2$ would be subtracted twice. To compensate we need to add back in $[x/p_1 \cdot p_2]$. The same would hold for any pair $p_i, p_j$. However, adding back $[x/p_i \cdot p_j]$ for each pair would cause another double counting effect, but in the opposite direction, when we consider the any 3 primes, say, $p_1, p_2$, and $p_3$. When we make the compensations for double counting the pairs $p_1 \cdot p_2, p_1 \cdot p_3$, and $p_2 \cdot p_3$, any multiples of $p_1 \cdot p_2 \cdot p_3$ would be added in twice, so we need compensate by subtracting $[x/p_1 \cdot p_2 \cdot p_3]$--and so on as indicted in formula (2).

While the theorem purports to be about counting *n*-primes, in reality it mostly counts composites involving the first *n* primes. Only when the numerous terms of the argument on the right side of the equation are computed and combined, and when that result is subtracted from $[x]$, is the number of *n*-primes that are ≤ x revealed. If the "$[x]$-" part of the formula is dropped, then all the remaining terms of the formula reveal the number of integers ≤ x that are divisible by one or more of the first *n* primes (in the case of formula (1)) or of the primes ε *M* (in the case of formula (3))

Legendre's formula is impractical after just a few iterations. The tenth iteration requires more than a thousand separate calculations and the 20[th] more than a million. It is a safe guess that Legendre himself never made any such calculations during his lifetime, long before today's computers were even imagined. But his formula was then, and still is today, valuable for the insights it provides about the patterns of the distribution of primes and composite integers. It allows the mind's eye to peer vast distances down the number line, not just to count primes but to better understand the patterns of their distribution.

**B. Introduction to Formulas for Counting Twins.** Our formulas for counting twins are analogous to Legendre's formula. However, they are necessarily more complex because they calculate the number of *pairs* of integers, *both* of which are relatively prime to the first *n* primes or to the primes in set *M*.

As with Legendre's formula, we situate our twin-counting functions in the domain of real numbers rather than integers, because intermediate stages of most of our calculation involve real numbers rather than integers. However, the results are always integers. Thus, our twin *n*-prime (and twin *M*-prime) counting functions map the set of real numbers, R, onto the set of positive integers Z, i.e.,

$$T_n(x): R \longrightarrow Z$$
$$T_M(x): R \longrightarrow Z$$
$$T\{p_a, p_b, p_c, \text{etc.}\}(x): R \longrightarrow Z$$
$$T\{p\}(x): R \longrightarrow Z$$

(4) $T_n(x)$ = the number of twin *n*-prime centers ≤ x, for the first *n* primes
(5) $T_M(x)$ = the number of twin *M*-prime centers ≤ x, for any set *M* of primes.
(6) $T\{p_a, p_b, p_c, \text{etc.}\}(x)$ = the number of twin *M*-prime centers ≤ x, for the primes specifically named. Thus, formulas (4) and (5) are identical when $M = \{p_a, p_b, p_c, \text{etc.}\}$; and
(7) $T\{p\}(x)$ = the number of twin *M*-prime centers ≤ x, for the set *M* consisting of the single prime, p.

These functions will count the number of integers, x, that are centered between two integers, (x – 1) and (x + 1), where neither of them is divisible by one of the first *n* primes, or by any of the primes of a given set *M*.

**C. A Formula for Counting Twins of a Single Prime.** We begin by establishing a formula for a single prime, p. It can be derived by reflecting on our first example, *p* = 7, illustrated in *Table 3, Centers of Twin {7}-primes*. There we see that the pattern of primes and composites is cyclical, of length 7. For each cycle, all integers are the centers of twin {7}-primes except the first and the second last, 1 and 6. This can readily be generalized to any prime, p. Thus, our fundamental formula for calculating the number of twin *n*-primes or *M*-primes ≤ x, for any prime *p* other than 2 and any real number x is:

(**8**) **T{p}(x) = [x] – [(x + 1)/p] – [(x + (p - 1))/p]**

For example, the formula for counting centers of twin pairs of integers both relatively prime to 5 would be:

(**9**) **T{5}(x) = [x] – [(x + 1)/5] – [(x + 4))/5]**

Similarly, the formula for counting centers of twin pairs of integers both relatively prime to 7 would be:

(**10**) **T{7}(x) = [x] – [(x + 1)/7] – [(x + 6))/7]**

NOTE: Theorem 1 does not hold for *p* = 2, since it would double count the number of integers to be subtracted from [x]. That is because for *p* = 2, inserting it as the value of x in formula (7) includes the two expressions [(x + 1)/2] and [(x + (2 - 1))/2] which are identical, and thus repetitive and double counting its value. Hence, for *p* = 2, the number of twin 2-primes is, more simply,

(**11**) **T{(2)}(x) = [x] – [(x + 1)/2]**

**D. Formulas for Counting Twins of Multiple Primes.** We now turn to the question of how many twin prime centers are ≤ x with respect to two or more primes. To illustrate the formula, we first consider the case of two primes, 5 and 7. On the surface, we may simply apply formula 1 twice, once for each of these two primes. So as a first shot, we might guess that

(**12**) **T{5, 7}(x) = [x] – ([(x+1)/5] – [(x+4)/5]) – ([(x+1)/7] – [(x+6)/7]).**

However, this would *not be correct*, because some values of T{5, 7}(x) would be implicated by both 5 and 7. This is illustrated in *Table 9, Centers of Twin {5}- and {7}-Primes*.

The first and fourth columns identify integers that are, respectively, centers of twin {5}-primes and twin {7}-primes. That is, for every number, x, in the first column, neither (x – 1) nor (x + 1) is divisible by 5. Similarly, for every number, x, in the fourth column, neither (x – 1) nor (x + 1) is divisible by 7.

However, the two middle columns highlight integers that are *not* centers of either twin {5}-primes or {7}-primes. For example, 16 is *not* the center of twin {5}-primes since (16 – 1) = 15 is divisible by 5; and 20 is *not* the center of twin {7}-primes, since (20 + 1) = 21 is divisible by 7.

**Table 9. Centers of Twin {5}- and {7}-Primes**

| Twin {5}-Prime Centers | Not Twin {5}-Prime Centers | Not Twin {7}-Prime Centers | Twin {7}-Prime Centers |
|---|---|---|---|
| (0) | | | |
| | 1 | 1 | |
| 2 | | | 2 |
| 3 | | | 3 |
| | 4 | | 4 |
| 5 | | | 5 |
| | 6 | 6 | |
| 7 | | | 7 |
| 8 | | 8 | |
| | 9 | | 9 |
| 10 | | | 10 |
| | 11 | | 11 |
| 12 | | | 12 |
| 13 | | 13 | |
| | 14 | | 14 |
| 15 | | 15 | |
| | 16 | | 16 |
| 17 | | | 17 |
| 18 | | | 18 |
| | 19 | | 19 |
| 20 | | 20 | |
| | 21 | | 21 |
| 22 | | 22 | |
| 23 | | | 23 |
| | 24 | | 24 |
| 25 | | | 25 |
| | 26 | | 26 |
| 27 | | 27 | |
| 28 | | | 28 |
| | 29 | 29 | |
| 30 | | | 30 |
| | 31 | | 31 |
| 32 | | | 32 |
| 33 | | | 33 |
| | 34 | 34 | |
| 35 | | | 35 |

Furthermore, as highlighted in the table, in four cases, namely when x = 1, 6, 29, and 34, *both* primes 5 and 7 are implicated—i.e., those four numbers are not centers of twin *M*-primes for either 5 or 7. For example, 29 is not the center of twin 7-primes, since (29 – 1) = 28 is divisible by 7; and it is also not the center of twin 5-primes, since (29 + 1) = 30 is divided by 5. Similarly, 1, 6, and 34 are bounded on one side or the other by multiples of both 5 and 7. Therefore equation (7) undercounts the value of $T_{\{5, 7\}}(x)$ for each of those four values by subtracting them twice. To correct for this undercount, we need to add back one repetition of each of the double counted values during every cycle of 35 integers. Thus, we need to substitute for equation (7) the following more elaborate formula described in Table 10.

A similar, but increasingly complex, adjustment to the formula is needed for each additional prime we introduce. For example, consider $T_{\{5, 7, 11\}}(x)$, the formula for the number of twin prime centers ≤ x whose wings, (x + 1) and (x – 1) are relatively prime to 5, 7, and 11.

Like Legendre's formula, the introduction of the third prime, 11, brings new complexities. It gives rise to two new prime pairs, for a total of three such pairs, whose products are needed to compensate for the double deletion of potential twin prime centers. However, also like Legendre's formula, this compensation in turn gives rise to the need of a triplet whose product needs to be added back in to restore twin prime centers that were deleted twice in the previous adjustment. And so on. By way of an example, the components of such a formula are identified in Table 13.

**Table 10. Number of Integers a ≤ x
That Are Centers of Integer Pairs (a – 1) and (a + 1)
Both Relatively Prime to 5 and 7**

| Prime Subsets | T{5, 7}(x) Twin Centers ≤ x |
|---|---|
| | [x] |
| {5} | – [(x+1)/5] – [(x+4)/5] |
| {7} | – [(x+1)/7] – [(x+6)/7] |
| {5, 7} | + [(x+1)/35] + [(x+6)/35] |
| | + [(x+29)/35] + [(x+34)/35] |

From these two examples, it is apparent that the formulas for counting *twin n*-primes, or *twin M*-primes, are even more complex than the comparable formulas for counting *n*-primes or *M*-primes. *Section C. A Formula for Counting Twins of a Single Prime*, hinted at this complexity by demonstrating that for a single prime, *p*, during each interval of consecutive *p* integers, two of them would be discarded as centers of twin {*p*}-primes, namely, [(x + 1)/p] and [(x + (p - 1))/p]. This complexity further unfolds as additional primes are taken under consideration, as reflected in *Table 11. Formula for the Number of Integers a ≤ x That Are Centers of Integer Pairs (a - 1) and (a + 1), Both Relatively Prime to 5, 7, and 11*.

**Table 11. Number of Integers a ≤ x
That Are Centers of Integer Pairs (a - 1) and (a + 1),
Both Relatively Prime to 5, 7, and 11**

| | | T{5, 7, 11}(x) | |
|---|---|---|---|
| Prime Subsets | Cycle Length | "Discarded Centers" (Integers Discarded as Possible Twin Centers Within Each Cycle) | Number of Twin Centers ≤ x |
| | | | [x] |
| {5} | 5 | 1, 4 | – [(x+1)/5] – [(x+4)/5] |
| {7} | 7 | 1, 6 | – [(x+1)/7] – [(x+6)/7] |
| {11} | 11 | 1, 10 | – [(x+1)/11] -- [(x+10)/11] |
| {5, 7} | 35 | 1, 6 | + [(x+1)/35] + [(x+6)/35] |
| | | 29, 34 | + [(x+29)/35] + [(x+34)/35] |
| {5, 11} | 55 | 1, 21 | + [(x+1)/55] + [(x+21)/55] |
| | | 34, 54 | + [(x+34)/55] + [(x+54)/55] |
| {7, 11} | 77 | 1, 34 | + [(x+1)/77] + [(x+34)/77] |
| | | 43, 76 | + [(x+43)/77] + [(x+76)/77] |
| {5, 7, 11} | 385 | 1, 34 | -- [(x+1)/385] -- [(x+34)/385] |
| | | 76, 111 | -- [(x+76)/385] -- [(x+111)/385] |
| | | 274, 309 | -- [(x+274)/385] -- [(x+309)/385] |
| | | 351, 384 | -- [(x+351)/385] -- [(x+384)/385] |

*Table 9, Centers of Twin {5}- and {7}-Primes,* along with its explanatory text, illustrated how certain integers, x, were discarded as potential centers of twin primes because the integers that immediately precede or follow them [e.g., (x + 1) or (x + 4); and (x + 1) or (x + 6)] could not serve as centers of twin *n*-primes because one or both members of each pair were divisible by either 5 or 7, or both. *Table 12, Integers Discarded as Twin Centers with Respect to 5, 7, and 11*, and related explanatory text, expand on this idea, highlighting the "discarded centers," that cannot be twin *M*-prime centers because one or both wings are divided by primes in subsets $M_{\{5,7\}}$, $M_{\{5,11\}}$, $M_{\{7,11\}}$, and $M_{\{5,7,11\}}$ of the set $M = \{5, 7, 11\}$.

As indicated in the first section, the numbers 1, 6, 29, and 34 are eliminated as twin *M*-prime centers for subset $M_{\{5, 7\}}$ of *M* because for each one, either its predecessor or successor, or both, are divided by either 5 or 7, or both. Analogously, 1, 21, 34, and 54 are eliminated as twin *M*-prime centers for the subset $M_{\{5, 11\}}$, as are 1, 34, 43, and 76 for subset $M_{\{7, 11\}}$.

The table also indicates that for each such pair, the respective patterns are repeated for intervals of lengths equal to their products, i.e., of lengths 5·7 = 35; 5·11 = 55; and 7·11 = 77. Thus, for example, for the pair M = {5, 7}, for any integer k, (35·k + 1), (35·k + 6), (35·k + 29), and (35·k + 34) are eliminated as centers of twin *M*-primes.

**Table 12. Integers Discarded as Twin Centers with Respect to Subsets of 5, 7, and 11**

| *M* = {5,7} | | *M* = {5, 11} | | *M* = {7, 11} | | *M* = {5, 7, 11} | |
|---|---|---|---|---|---|---|---|
| Integers discarded as twin centers and associated twin pairs | One or both of pair divided by 5 or 7 | Integers discarded as twin centers and associated twin pairs | One or both of pair divided by 5 or 11 | Integers discarded as twin centers and associated twin pairs | One or both of pair divided by 7 or 11 | Integers discarded as twin centers and associated twin pairs | One or both of pair divided by 5, 7, or 11 |
| 1 | 0 / 5, 7 | 1 | 0 / 5, 11 | 1 | 0 / 7, 11 | 1 | 0 / 5, 7, 11 |
|  | 2 |  | 2 |  | 2 |  | 2 |
| 6 | 5 / 5 | 21 | 20 / 5 | 34 | 33 / 11 | 34 | 33 / 11 |
|  | 7 / 7 |  | 22 / 11 |  | 35 / 7 |  | 35 / 5, 7 |
| 29 | 28 / 7 | 34 | 33 / 11 | 43 | 42 / 7 | 76 | 75 / 5 |
|  | 30 / 5 |  | 35 / 5 |  | 44 / 11 |  | 77 / 7, 11 |
| 34 | 33 | 54 | 53 | 76 | 75 | 111 | 110 / 5, 11 |
|  | 35 / 5,7 |  | 55 / 5, 11 |  | 77 / 7, 11 |  | 112 / 7 |
| Cycle repeats every 35 integers | | Cycle repeats every 55 integers | | Cycle repeats every 77 integers | | 274 | 273 / 7 |
| | | | | | | | 275 / 5, 11 |
| | | | | | | 309 | 308 / 7, 11 |
| | | | | | | | 310 / 5 |
| | | | | | | 351 | 350 / 5, 7 |
| | | | | | | | 352 / 11 |
| | | | | | | 384 | 383 |
| | | | | | | | 385 / 5, 7, 11 |
| | | | | | | Cycle repeats every 385 integers | |

The last section of the table illustrates similar results for the set *M* = {5, 7, 11}. However, because three primes are involved rather than two, there are eight integers discarded as centers of twin *M*-primes with respect to those three primes. In this case, the cycle repeats itself every 385 integers, the product of 5, 7, and 11.

**E. A General Formula for Counting Pairs of Twins.** We are now in a position to describe a formula for counting pairs of twin primes, twin *n*-primes, and twin *M*-primes that is analogous to Legendre's formula for counting all primes.

**Definition: Discarded Twin Centers.** We begin with sets that we call "discarded twin centers." These are illustrated in column 3 of Table 11 and are more fully explained in connection with Table 12. We now define them more specifically and illustrate how they are used in counting pairs of twin *n*-primes or twin *M*-primes.

Let *M* be a finite set of primes. Define:
"$D_M$": the set of "discarded *M*-prime centers," i.e., integers x such that
x ε $P_M$ (i.e., $0 \leq x \leq \Pi p_i$, $p_i$ ε M) and that each member of M divides either (x – 1) or (x + 1).

Analogous definitions and notations are:
"Dn": the set of discarded *M*-prime centers for the first *n* primes. We call these "discarded *n*-prime centers."
"$D\{p_a, p_b, p_c \ldots\}$": the set of discarded *M*-prime centers for the set M = {$p_a, p_b, p_c \ldots$}.
"$D\{p_i\}$": the set of discarded *M*-prime centers for the set M consisting of the single prime, $p_i$.

Examples of discarded *M*-prime centers are displayed in the third column of *Table 11. Formula for the Number of Integers a ≤x That Are Centers of Integer Pairs (a - 1) and (a + 1), Both Relatively Prime to 5, 7, and 11*, and in the subheadings entitled "Discarded centers and pairs" of *Table 12. Integers Discarded as Twin Centers with Respect to Subsets of 5, 7, and 11.*

**Function: Number of Discarded Twin Centers Relative to a Set M.** For a given set M of primes, define the function $D_M(x)$, mapping positive real numbers R onto the non-negative integers Z

$$D_M(x): R \longrightarrow Z$$

as the number of integers less than or equal to x that are discarded as twin centers relative to the set M. Thus,

(13)  $D_M(x) = \Sigma[(x + z)/\Pi p_i]$ for all $p_i$ ε M and z ε $D_M$

For example, from *Table 11. Number of Integers a ≤x that Are Centers of Integer Pairs (a - 1) and (a + 1), Both Relatively Prime to 5, 7, and 11* we see that

$D\{5\}(x) = [(x+1)/5] + [(x+4)/5]$

$D\{5, 7\}(x) = [(x+1)/35] + [(x+6)/35] + [(x+29)/35] + [(x+34)/35]$

$D\{5, 7, 11\}(x) = [(x+1)/385] + [(x+34)/385] + [(x+76)/385] + [(x+111)/385] +$
    $[(x+274)/385] + [(x+309)/385] + [(x+351)/385] + [(x+384)/385]$

**Twin Counting Functions.** With this definition and examples in mind we now define our function Tn(x) for counting twin centers ≤ x as follows.

(14)  $Tn(x) = [x] – \Sigma(D\{p_i\}(x)) + \Sigma(D\{p_i, p_j\}(x)) - \Sigma(D\{p_i, p_j, p_k\}(x)) + \cdots$
    $\cdots \pm \Sigma(D\{p_1, p_2, \ldots p_n\}(x))]$
sums taken over all the singles, then all the doubles, triples, etc. of the first n primes and alternately subtracted from and added to [x].

For any integer *n* and any real number x, this formula returns the number of integers ≤ x that are centers of twin n-primes, with respect to the first *n* primes.

We may generalize this formula to integers relatively prime to all primes within a given set *M* of primes, i.e.

(15)  $T_M(x) = [x] - \Sigma(D_{\{p_i\}}(x)) + \Sigma(D_{\{p_i, p_j\}}(x)) - \Sigma(D_{\{p_i, p_j, p_k\}}(x)) + \cdots$
$\cdots \pm \Sigma(D_{\{p_i \cdot p_j \cdot p_k \cdots\}}(x))]$

the sums taken over each, then all the doubles, triples, etc. of the primes ε M.

The formats of these two formulas are intended to mimic Legendre's formula for counting primes[viii]. We may express the same result more efficiently using function (13), $D_M(x)$, for discarded twin centers as follows:

Let M be a finite set of primes; and let $m_1, m_2, m_3, \ldots$ etc. be any enumerated list of all the subsets of *M*. Then

(16)  $T_M(x) = [x] - \Sigma((-1)^k (D_{m_i}(x)))$
where k = 1 if the cardinality of $m_i$ is even, and 2 if odd.

As with Legendre's formula, while formulas (15) and (16) are valid and precise, they are impractical for making manual calculations after just a few iterations. By way of comparison, ten iterations of the twin counting formula require more than 59,000 calculations; and 20 iterations require almost 3.5 billion. The 59,000 are of course still a pittance of the 100 billion integers in the natural domain of the first 10 odd primes, and again, the pattern of the distribution of those twins repeats throughout the number line within endless cyclical multiples of these first 100 billion integers.)

**5. Meissel/Lehmer-inspired Twin Counting Functions**

**A. Meissel/Lehmer's Formula.** Here we will examine whether Meissel/Lehmer's formula for counting primes can be adapted to construct a new formula for counting *twin* primes, *twin n*-primes, and *twin M*-primes.

Meissel-Lehmer's first introduced by Ernst Meissel in 1870, is

(17)  $f_n(x) = f_{n-1}(x) - f_{n-1}(x/p_n)$
where $f_n(x)$ is the number of *n*-primes.

The principle underlying this formula is that an integer *x* that is relatively prime to the first *n - 1* primes is also relatively prime to the first *n* primes unless it is divisible by n[th] prime, $p_n$. But that happens if and only if it also has a divisor that is ≤ $x/p_n$.

Unfortunately, the situation is not analogous for *twin n*-primes. They lose their status as twins if *one* or *both* members of the pair is a composit integer, and even if a divisor of *either* member of the pair is *any* prime, not just a *twin* prime. Furtermore, we have cast our analysis of twin *n*-primes and twin M-primes in terms of their *centers*. Thus, our goal here is to count which, and how many, twin (*n-1*)-prime, or twin *M*-prime, *centers* survive as twin *centers* with the introduction of an additional prime. Adapting Meissel-Lehmer's prime counting formula to account for these features of twin *n*-primes or *M*-primes introduces so many complexities that the value of a single formula as elegant as Meissel-Lehmer's prime counting formula is lost.

**Inductive Method for Counting Twins.** Despite the lack of a straightforward adaptation of Meissel/Lehmer's formula for counting twin primes, we may gain insights about the effects of introducing an additional prime into the mix by examining *Table 15, Incremental Growth in the Number of Integer Pairs (a - 1) and (a + 1) Relatively Prime to 5, 7, and 11*. The content of this table is identical to that of *Table 13*, but the order and grouping of the rows have been moved to focus on the changes that occur with the incremental introduction of each of the three primes, 5, 7, and 11.

The first row reveals that two calculations ( $- [(x+1)/5] - [(x+4)/5]$ ) are needed to account for the number of integers discarded as potential twin centers by the single prime, 5. However, when a second prime, e.g., 7, is introduced, six additional calculations are required--two to account for the number of integers that are discarded by 7, and then four more to compensate for double counting the combined effects of 5 and 7. The introduction of a third prime, 11, introduces even more complications, first to account for the centers discarded by 11; then to compensate for double counting by the combined effects of 5 and 11 and of 7 and 11; but then again to account for double counting the effects of those two compensations.

These modest examples reveal that things get complicated very quickly. The number of calculations needed to compute the number of twin centers ≤ x grows exponentially by a factor of 3. For example, *Table 13, Incremental Growth in the Number of Integer Pairs (a - 1) and (a + 1) Relatively Prime to 5, 7, and 11,* reveals that, if we consider each bracketed sub-calculation in the fourth column (e.g., $[(x+1)/5]$ ) to be a *single* calculation, then 3 such calculations are needed to compute the number of twins **≤ x** with respect to the single prime, 5; 8 for the pair, 5 and 7; and 27 for the triplet, 5, 7, and 11. Along these same lines, 1,000 calculations would be needed to compute the number of twins **≤ x** with respect to any 10 primes; 8,000 calculations for 20 primes; 27,000 for 30 primes, etc.

**Table 13. Incremental Growth in the Number of Integer Pairs (a - 1) and (a + 1) Relatively Prime to 5, 7, and 11**

| | | $T_{\{5, 7, 11\}}(x)$ | |
|---|---|---|---|
| **Prime Subsets** | **Cycle Length** | **"Discarded Centers" (Integers Discarded as Possible Twin Centers Within Each Cycle)** | **Number of Twin Centers ≤ x** |
| | | | [x] |
| **{5}** | 5 | 1, 4 | $- [(x+1)/5] - [(x+4)/5]$ |
| **{7}** | 7 | 1, 6 | $- [(x+1)/7] - [(x+6)/7]$ |
| **{5, 7}** | 35 | 1, 6<br>29, 34 | $+ [(x+1)/35] + [(x+6)/35]$<br>$+ [(x+29)/35] + [(x+34)/35]$ |
| **{11}** | 11 | 1, 10 | $- [(x+1)/11] -- [(x+10)/11]$ |
| **{5, 11}** | 55 | 1, 21<br>34, 54 | $+ [(x+1)/55] + [(x+21)/55]$<br>$+ [(x+34)/55] + [(x+54)/55]$ |
| **{7, 11}** | 77 | 1, 34<br>43, 76 | $+ [(x+1)/77] + [(x+34)/77]$<br>$+ [(x+43)/77] + [(x+76)/77]$ |
| **{5, 7, 11}** | 385 | 1, 34<br>76, 111<br>274, 309<br>351, 384 | $-- [(x+1)/385] -- [(x+34)/385]$<br>$-- [(x+76)/385] -- [(x+111)/385]$<br>$-- [(x+274)/385] -- [(x+309)/385]$<br>$-- [(x+351)/385] -- [(x+384)/385]$ |

With such an exponential growth rate, the number of calculations needed to compute $T_M(x)$ would quickly exceed the capacity of commonly available computers. Of course, the field of mathematics has in recent years opened up to and kept pace with the development of super computers that are now widely used in counting prime numbers. However, the purpose of this paper is not to promote or facilitate such prime-counting, but rather to shed light on the theoretical underpinnings of the distributions of twin primes, *n*-primes, and *M*-primes within the remnants of the sieve of Eratosthenes.

**B. Calculating the Discarded Centers.** The central focus of our formulas for counting the number of n-primes and M-primes and in describing their distribution along the number line are the "discarded centers." For a set of primes, *M*, we have defined the discarded centers with respect to *M* as integers z such that $0 \leq z \leq \Pi p_i$, (pi ε *M*), and that each member of *M* divides either (z – 1) or (z + 1). Thus, the pair (z - 1), (z + 1) cannot be twin M-primes, since one or both members of the pair are divisible by at least one $p_i$.

For example, in Tables 12 and 13, consider the subset M = {5, 11}. As both tables indicate, its "discarded centers" are 1, 21, 34, and 54. This means that they cannot be centers of twin primes because 5 and 11 divide one or both integers above or below them. For example,

(18)  (1 – 1) = 0, is divisible by both 5 and 11
(1 + 1) = 2

(21 - 1) = 20, is divisible by 5,
(21 + 1) = 22, is divisible by 11

(34 – 1) = 33, is divisible by 11
(34 + 1) = 35, is divisible by 5

(54 – 1) = 53
(54 + 1) = 55, is divisible by both 5 and 11

Thus, none of these "discarded centers" can be the centers of twin *M*-primes where M = {5, 11}. This pattern repeats itself at intervals of 5 x 11 = 55 across the entire number line. Analogous "discarded centers" and associated cycles appear in tables 12 and 13 for all the subsets of the set M = {5, 7, 11}.

These centers play a key role in counting twin-*M* primes that is similar to what the prime numbers and the various sets and products of them play in Legendre's formula for counting twin primes. (See section 3)

At this point, an open question is: how were these "discarded centers" determined? A general answer is that they were calculated using Diophantine equations of two unknowns, which typically involve integers in their formulation and solutions. A full discussion of Diophantine questions is beyond the scope of this paper. However, the following examples may give at least a taste of how such equations can be used to determine the discarded centers of twin n-primes or M-primes.

One method is to identify the twin centers associated with the set $M = \{5, 11\}$ and then see which of them survive as twin centers associated with the larger set $M = \{5, 7, 11\}$.

For example, consider 34, which, as noted above, is a "discarded center" with respect to the set M = {5, 11}. This is due to the fact each member of the set M divides one of its wings. For example, 7 divides ((34 + 1) = 35) and 11 divides ((34 - 1) = 33). It is also a discarded center with respect to the set M = {5, 11} within *each successive multiple of 55*, as indicated in *Table 12. Integers Discarded as Twin Centers with Respect to Subsets of 5, 7, and 11.* This cyclical result, for 7 cycles, is displayed in *Table 14. Example of Derivation of Twin M-Prime Center From M-{5, 11} to M-{5, 7, 11}*.

**Table 14. Example of Derivation of Twin M-Prime Center From M-{5, 11} to M-{5, 7, 11}**

| 7 Cycles of length 55 | x | x + 55 | x + 110 | x + 165 | x + 220 | x + 275 | x + 330 |
|---|---|---|---|---|---|---|---|
| Discarded Center | 34 | 89 | 144 | 199 | 254 | 309 | 364 |
| Wings | 33 35 | 88 90 | 143 145 | 198 200 | 253 255 | 308 310 | 363 365 |
| Divisors of Wings by 5 or 11 | 11 5 | 11 5 | 11 5 | 11 5 | 11 5 | 11 5 | 11 5 |
| Division of Wings by 7 | 4.7 **5** | 12.6 12.9 | 20.4 20.7 | 28.3 28.6 | 36.1 36.4 | **44** 44.3 | 51.9 52.1 |

We may note that in the last row, in only two of the cycles, the first and the 6[st], does the division of the wings by 7 yield an integer - 5 in the first column and 44 in the sixth[th]. These correspond to two of the discarded centers, 34 and 309, identified in *Table 12. Integers Discarded as Twin Centers with Respect to Subsets of 5, 7, and 11.*

Repeating this process for each of the three remaining discarded centers of M = {5, 11}, namely 1, 21, and 54, will yield the remaining 6 discarded centers of M = {5, 7, 11}.

### 6. Euler-like ϕ-function for Twins

For any integer $x = p_1^{k1} p_2^{k2} p_3^{k3} \cdots p_n^{kn}$, Euler's ϕ function describes the number of integers ≤ x that are relatively prime to it. The formula is written in several equivalent formats, including:

18. $\phi(x) = x \cdot \dfrac{(p_1 - 1)(p_2 - 1)(p_3 - 1) \cdots (p_n - 1)}{p_1 \cdot p_2 \cdot p_3 \cdots p_n}$

where $p_1, p_2, \ldots p_n$ are the set M of unique prime dividers of x.

Here we wish to introduce a conceptually parallel function, T(x), which returns the number of integers $a \leq x$ that are *centers* of twin integers (a - 1) and (a + 1), with both relatively prime to x. Like ϕ (x), T(x) may be written in several equivalent formats, including:

19. $T(x) = x \cdot (p_1 - 2)(p_2 - 2)(p_3 - 2) \cdots (p_n - 2)$

$$p_1 \cdot p_2 \cdot p_3 \cdots p_n$$
where $p_1, p_2, \ldots p_n$ are the set of unique prime dividers of x.

This result follows directly from Theorem 4a (in the background section) with respect to the set M consisting of those unique prime dividers of x.

**7. Conclusion**

Broadly speaking, the distributions of twin n-primes and twin *M*-primes follow patterns remarkably similar to that of *all* n-primes and M-primes, and to some extent that of twin primes themselves. Still, our vision and grasp of those infinite worlds is still murky, especially that of the twin primes, whether finite or infinite.

# Appendix

**A. Key Findings from Cycles and Patterns in the Sieve of Eratosthenes**
G. Grob and M. Schmitt; 28 Oct. 2019 (https://arxiv.org/abs/1905.03117)

The above paper introduces the concept of *n*-primes and *M*-primes and thereby provides a starting point and framework for our analysis of *twin n*-primes and *twin M*-primes. The following theorems from that paper are most germane to our analysis of the twins.

**For *n*-primes (integers not divisible by any of the first *n* primes)**

- **Theorem 1.** The distribution of *n*-primes is repeated every interval of length $\Pi p_i$, $1 \leq i \leq n$. That is, if x is an *n*-prime, so is $K \cdot \Pi p_i + x$ for every integer K and for $p_i$, $1 \leq i \leq n$.
- **Theorem 2.** N-primes are distributed symmetrically within every such interval $\Pi p_i$. That is, if x is an *n*-prime, so is $\Pi p_i - x$. By Theorem 1, this also holds true for every cycle $K \cdot \Pi p_i$ for every integer K.
- **Theorem 3.** The number of *n*-primes in every $\Pi p_i$ interval = $\Pi(p_i - 1)$, for $1 \leq i \leq n$
- **Corollary 1.** For any integer $K \geq 0$, $f_n(K \cdot \Pi p_i \pm x) = f_n(K \cdot \Pi p_i) \pm f_n(x) = K \cdot \Pi(p_i - 1) \pm f_n(x)$, for $1 \leq i \leq n$
- **Theorem 4.** Twin n-primes are pairs consisting of two *n*-primes, (a – 1) and (a +1). The number of such pairs $\leq \Pi p_i = \Pi(p_i - 2)$, $1 \leq i \leq n$ ($p_i \neq 2$).
- **Corollary 2.** For any integer $a \neq 0$ (mod $p_i$, $1 \leq i \leq n$) and any integer $K \geq 0$ there are $K \cdot \Pi(p_i - 2)$, $1 \leq i \leq n$, ($p_i \neq 2$) integers $x \leq K \cdot \Pi p_i$ such that $(x - a)$ and $(x + a)$ are both *n*-primes.
- **Theorem 8, Meissel's formula.** (Repeated here for ease of reference) $f_n(x) = f_{n-1}(x) - f_{n-1}(x/p_i)$ for any $p_i$, $1 \leq i \leq n$.

**For *M*-Primes (integers not divisible by any prime in the set M of primes)**

- **Theorem 1a.** The distribution of *M*-primes is repeated every $\Pi p_i$ interval, $p_i \varepsilon M$. That is, if x is an *M*-prime, so is $K \cdot \Pi p_i + x$ for every integer K and for $p_i \varepsilon M$
- **Theorem 2a.** *M*-primes are distributed symmetrically within $P_M$. That is, if x is an *M*-prime, so is $\Pi p_i - x$. By Theorem 1, this also holds true for every cycle $K \cdot \Pi p_i$ for every integer K.
- **Theorem 3a.** The number of *M*-primes in $P_M = \Pi(p_i-1)$, $p_i \varepsilon M$
- **Corollary 1a.** For any integer $K \geq 0$, $f_M(K \cdot \Pi p_i \pm x) = f_M(K \cdot \Pi p_i) \pm f_M(x) = K \cdot \Pi(p_i - 1) \pm f_M(x)$, $p_i \varepsilon M$
- **Theorem 4a.** Twin M-primes are pairs M-primes, (a – 1) and (a +1). The number of such pairs $\leq \Pi p_i = \Pi(p_i - 2)$, ($p_i \varepsilon M$), ($p_i \neq 2$).
- **Corollary 2a.** For any integer $a \neq 0$ (mod $p_i$, $p_i \varepsilon M$) and any integer $K \geq 0$ there are $K \cdot \Pi(p_i - 2)$, $p_i \varepsilon M$, ($i \neq 2$) integers $x \leq K \cdot P_n$ such that $x - a$ and $x + a$ are both *M*-primes.
- **Theorem 8a, Meissel's formula generalized.** (Repeated here for ease of reference) $f_M(x) = f_{M-pi}(x) - f_{M-pi}(x/p_i)$ for any $p_i \varepsilon M$. (Note, Meissel's formula is traditionally expressed in terms of the first *n* primes. Here, we generalizes it to arbitrary set *M* of primes, including infinite sets.

**B. Notations.** Following are the notations use throughout this paper.

1) $\pi(x)$ = the number of primes $\leq x$
2) $f_n(x)$ = the number of $n$-primes $\leq x$
3) $f_M(x)$ = the number of integers $\leq x$ that are relatively prime to all $p_i \, \varepsilon \, M$, where M is a finite or infinitely countable set of primes
4) $f\{p_a, p_b, p_c \ldots\}(x)$ = the number of integers $\leq x$ that are relatively prime to all the specifically designated primes, $p_i$
5) $f\{p_i\}(x)$ = the number of integers that are $\leq x$ and relatively prime to the single prime, $p_i$. (This is a special case of (4))
6) $T(x)$ = the number of twin prime centers $\leq x$.
7) $T_n(x)$ = the number of twin $n$-primes centers $\leq x$
8) $T_M(x)$ = the number of twin $M$-prime centers $\leq x$
9) $T\{p_a, p_b, p_c \ldots\}(x)$ = the number of twin $M$-prime centers $\leq x$ for a specific set $M$ of primes = $\{p_a, p_b, p_c \ldots\}$, where M is a finite or infinitely countable set of primes
10) $T\{p_i\}(x)$ = the number of twin prime centers that are $\leq x$ and relatively prime to the single prime, $p_i$. (This is a special case of (9))
11) Define $P_n$ as the set of integers $0 \leq x \leq \Pi p_i$, $i \leq n$
12) Define $P_M$ as the set of integers $0 \leq x \leq \Pi p_i$, $p_i \, \varepsilon \, M$.
13) Define $D_M$ as the set of "discarded $M$-prime centers," i.e., integers x such that $x \, \varepsilon \, P_M$ and that each member of $M$ divides either $(x - 1)$ or $(x + 1)$.
14) Define $D_n$ as the set of discarded $M$-prime centers for the first $n$ primes. We call these "discarded $n$-prime centers."
15) Define $D\{p_a, p_b, p_c \ldots\}$ as the set of discarded $M$-prime centers for the set $M = \{p_a, p_b, p_c \ldots\}$.
16) Define $D\{p_i\}$ as the set of discarded $M$-prime centers for the set M consisting of the single prime, $p_i$.

---

[i] G. Grob and M. Schmitt, Cycles and Patterns in the Sieve of Eratosthenes, arXiv: 1905.03117[math.GM] (revised 10-28-19), available at https://arxiv.org/abs/1905.03117

[ii] Adrien-Marie Legendre's (1752 – 1833) formula, using the inclusion-exclusion principle and adapted to n-primes, is
$f_n(x) = [x] - \Sigma [x/p_i] + \Sigma [x/p_i \cdot p_j] - \Sigma [x/p_i \cdot p_j \cdot p_k] + \ldots \Sigma [x/p_1 \cdot p_2 \cdot p_3 \ldots p_n]$ where $f_n(x)$ is the number of n-primes $\leq x$, and $p_i$ is the $i^{th}$ prime, $1 \leq i \leq n$.

[iii] Meissel-Lehmer's formula, first formulated by Ernst Meissel in 1870, is
$f_n(x) = f_{n-1}(x) - f_{n-1}(x/p_n)$, where $f_n(x)$ is the number of n-primes.

[iv] Euler's $\phi$-function, $\phi(n)$, is the number of non-negative integers that are less than, but relatively prime to, n.

[v] Ibid

[vi] Numerous texts describe the inclusion-exclusion principle. It is used to count the number of elements that satisfy at least one of several properties while avoiding double counting elements that satisfy more than one property. See, for example the explanation in Wikipedia.